\documentclass[12pt,letterpaper,twoside]{article}

\usepackage{amsmath,amssymb}
\usepackage{amsthm}
\usepackage{stmaryrd}
\usepackage{fancyhdr}
\usepackage{times}
\usepackage{mathrsfs} 
\usepackage{sectsty}
\usepackage[pdftex,dvips]{geometry}

\input diagxy.tex

\sectionfont{\normalsize\bfseries\scshape\centering\MakeUppercase }

\subsectionfont{\normalsize\bfseries\scshape }

\makeatletter
  \def\@seccntformat#1{\csname the#1\endcsname.\quad}
\makeatother

\newtheoremstyle{fact}
     {\topsep}
     {\topsep}
     {\slshape}
     {}
     {\bfseries}
     {}
     { }
     {\thmname{#1}\thmnumber{ #2.}\thmnote{ \rm (#3)}}

\newtheorem{theorem}{Theorem}[section]
\newtheorem*{theorem*}{Theorem} 
\newtheorem{lemma}[theorem]{Lemma}
\newtheorem{proposition}[theorem]{Proposition}
\newtheorem{corollary}[theorem]{Corollary}
\newtheorem{problem}{Problem}

\theoremstyle{definition}

\newtheorem{remark}[theorem]{Remark}
\newtheorem*{remark*}{Remark}

\newtheorem*{question*}{Question}
\newtheorem*{examples*}{Examples}  
\newtheorem{example}[theorem]{Example}
\newtheorem*{example*}{Example}

\theoremstyle{fact}

\newtheorem{ftheorem}[theorem]{Theorem}

\def\proofont{\fontseries{bx}\fontshape{sc}\selectfont}
\def\proofname{Proof. }

\makeatletter
\renewenvironment{proof}[1][\proofname]{\par
  \normalfont
  \topsep6\p@\@plus6\p@ \trivlist
  \item[\hskip\labelsep\noindent\proofont #1]\ignorespaces
}{%
  \qed\endtrivlist
}
\makeatother

\author{G\'abor Luk\'acs
\thanks{I gratefully acknowledge the generous financial support received
from the Alexander von Humboldt Foundation, the Killam Trusts, and 
Dalhousie University that enabled me to do this research.}}
   
\title{Almost maximally almost-periodic group topologies determined
by $T$-sequences}

\begin{document}

\makeatletter
\let\mytitle\@title
\chead{\small\itshape G. Luk\'acs / \mytitle }
\fancyhead[RO,LE]{\small \thepage}
\makeatother

\maketitle

\def\thanks#1{} 

\thispagestyle{empty}

\begin{abstract}
A sequence $\{a_n\}$ in a group $G$ is a {\em $T$-sequence} if there is a 
Hausdorff group topology $\tau$ on $G$ such that 
$a_n\stackrel\tau\longrightarrow 0$. In this paper, we provide several
sufficient conditions for a sequence in an abelian group to be 
a $T$-sequence, and investigate special sequences in 
the Pr\"ufer groups $\mathbb{Z}(p^\infty)$.
We show that for $p\neq 2$, there is a Hausdorff group topology
$\tau$ on $\mathbb{Z}(p^\infty)$ that is determined by a $T$-sequence,
which is close to being maximally almost-periodic---in other words, 
the von Neumann radical
$\mathbf{n}(\mathbb{Z}(p^\infty),\tau)$ is a non-trivial finite subgroup.
In particular, $\mathbf{n}(\mathbf{n}(\mathbb{Z}(p^\infty),\tau))
\subsetneq \mathbf{n}(\mathbb{Z}(p^\infty),\tau)$.
We also prove that the direct sum of any infinite family of finite 
abelian groups admits a group topology determined by a $T$-sequence with 
non-trivial finite von Neumann radical.
\end{abstract}

\section{Introduction}

Given a set $X$, a point $x_0 \in X$, and a sequence $\{x_n\}$ of
distinct elements in $X$, one can easily construct 
a Hausdorff topology $\tau$ on $X$ such that 
$x_n\stackrel \tau \longrightarrow x_0$. This is, however, not the 
case for groups and group topologies, as Example~\ref{ex:p-n:non-T}
below demonstrates.

Following Zelenyuk and Protasov \cite{PZMono} \& \cite{ZelProt}, 
who were the first to investigate this type of question, we say that 
a sequence $\{a_n\}$ in a group $G$ is a {\em $T$-sequence} if there is a 
Hausdorff group topology $\tau$ on $G$ such that 
$a_n\stackrel \tau \longrightarrow 0$. In this case,
the group $G$ equipped with the finest group topology with this property
is denoted by $G\{a_n\}$. A similar notion exists for filters, 
in which case one speaks of a {\em $T$-filter}.
Zelenyuk and Protasov characterized $T$-sequences and $T$-filters  in
abelian groups \cite[2.1.3, 2.1.4]{PZMono} \& \cite[Thm.~1, Thm.~2]{ZelProt}, 
and studied the topological  properties of $G\{a_n\}$, where $\{a_n\}$ 
is a $T$-sequence. (They also present a characterization of
$T$-filters in non-abelian groups in \cite[3.1.4]{PZMono}.)

\begin{example} \label{ex:p-n}
For a prime number $p$,  let $A=\mathbb{Z}(p^\infty)$ be the Pr\"ufer
group. It can be seen as the subgroup of $\mathbb{Q}/\mathbb{Z}$ generated
by the elements of $p$-power order, or the group formed by all
$p^n$th roots of unity in $\mathbb{C}$. For $e_n= \frac 1 {p^n}$, $\{e_n\}$ is 
clearly a  $T$-sequence in $A$, because $e_n \longrightarrow 0$
in the subgroup topology that $A$ inherits
from $\mathbb{Q}/\mathbb{Z}$.
\end{example}

\begin{example} \label{ex:p-n:non-T}
Keeping the notations of Example~\ref{ex:p-n}, set 
$a_n = -\frac 1 p + \frac 1 {p^n} = -e_1 +e_n$. If 
$a_n \stackrel \tau \longrightarrow 0$ for some group topology $\tau$, then 
also $e_{n-1} = p a_n \stackrel \tau \longrightarrow 0$, and therefore,
$e_1 = p a_{n+1} - a_n \stackrel \tau \longrightarrow 0$. Hence,
$\tau$ cannot be Hausdorff, and so $\{a_n\}$ is not a $T$-sequence in $A$.
\end{example}

Every topological group $G$
admits a ``largest" compact Hausdorff group $bG$ and a continuous
homomorphism $\rho_G \colon G \rightarrow bG$ such that
every continuous homomorphism $\varphi\colon G\rightarrow K$ 
into a compact Hausdorff group $K$ factors uniquely through 
$\rho_G$:
\begin{equation}
\bfig
\Vtriangle(0,0)/->`->`<--/<300,350>%
[G`K`bG;\varphi`\rho_G`\exists!\tilde\varphi]
\efig
\end{equation}
The group $bG$ is called the {\em Bohr-compactification} of $G$.
The image $\rho_G(G)$ is dense in $bG$. The kernel of
$\rho_G$ is called the {\em von Neumann radical} of $G$, and
denoted by $\mathbf{n}(G)$. One says that $G$ is
{\em maximally almost-periodic} if $\mathbf{n}(G)=1$, and
{\em minimally almost-periodic} if $\mathbf{n}(G)=G$ (cf.~\cite{NeuWig}). 
For an abelian topological group $A$, 
let $\hat A = \mathscr{H}(A,\mathbb{T})$ be the Pontryagin dual of 
$A$---in other words, the group of {\em continuous characters} 
of $A$ (i.e., continuous homomorphisms $\chi\colon A\rightarrow\mathbb{T}$, 
where $\mathbb{T}=\mathbb{R}/\mathbb{Z}$), equipped with the compact-open
topology. It follows from the famous Peter-Weyl Theorem
(\cite[Thm.~33]{Pontr}) that the Bohr-compactification
of $A$ can be quite easily computed:
$bA = \widehat{\hat A_d}$, where $\hat A_d$
stands for the group $\hat A$ with the discrete topology. Thus,
$\mathbf{n}(A)=\bigcap\limits_{\chi\in\hat A} \ker \chi$.

$T$-sequences turn out to be a very useful tool for constructing 
``pathological" examples. For example, Zelenyuk and Protasov used 
$T$-sequences to show (independently of Ajtai, Havas and Koml\'os 
\cite{AHK}) that every infinite abelian group admits a non-maximally 
almost-periodic Hausdorff group topology (cf.~\cite[2.6.4]{PZMono}, 
\cite[Thm.~16]{ZelProt}). There are 
plenty of examples of minimally (or maximally) almost-periodic groups
(cf.~\cite{NeuWig}, \cite{Nienhuys}). Nevertheless, it appears that no
example is known for a Hausdorff topological group $G$ whose 
von Neumann radical $\mathbf{n}(G)$ is non-trivial and finite.
We call such groups {\em almost maximally almost-periodic}.
This raises the following question:

\begin{problem} \label{prob:nG}
Which abelian groups $A$ admit a $T$-sequence $\{a_n\}$ such that
$A\{a_n\}$ is almost maximally almost-periodic?
\end{problem}

If $\mathbf{n}(G)$ is non-trivial and finite, then 
$\mathbf{n}(\mathbf{n}(G))=1$, and thus 
$\mathbf{n}(\mathbf{n}(G))\neq \mathbf{n}(G)$,
which leads to a second problem:

\begin{problem} \label{prob:nnG}
Which abelian groups $A$ admit a $T$-sequence $\{a_n\}$ such that
$\mathbf{n}(\mathbf{n}(A\{a_n\}))$ is strictly contained 
in $\mathbf{n}(A\{a_n\})$?
\end{problem}

Note that since $\mathbf{n}$ is productive, one may wish to focus on {\em 
algebraically} (respectively, {\em topologically}) {\em directly 
indecomposable groups}, that is, groups that cannot be expressed as an 
algebraic (respectively, topological) direct product of two of its proper 
subgroups.

\bigskip

Our ultimate goal in this paper is to present ample non-isomorphic 
algebraically directly indecomposable almost maximally almost-periodic
Hausdorff abelian groups (i.e., having non-trivial finite von Neumann 
radical), whose topology is determined by a $T$-sequence.
It will also show that Problems~\ref{prob:nG} and 
\ref{prob:nnG} are meaningful.

This aim is carried out according to the following structure:
In section~\ref{sect:abelian}, several results that provide sufficient 
conditions for a sequence in an abelian group to be a $T$-sequence are 
presented (Theorem~\ref{thm:T:ord}). In section~\ref{sect:sum}, 
a partial answer to Problem~\ref{prob:nG} is provided, namely, 
we prove that the direct sum of any infinite family of finite abelian groups 
admits an almost maximally almost-periodic group topology determined by a 
$T$-sequence (Theorem~\ref{thm:sum:main}). Groups of this form certainly 
fail to be algebraically directly indecomposable, and they need not be
be topologically directly indecomposable either. Thus, they fall short of 
our ultimate goal. In section~\ref{sect:Pruefer}, 
special sequences of Pr\"ufer groups are investigated 
(Theorem~\ref{thm:Pruefer:spec}), and we prove that 
for $p\neq 2$, $\mathbb{Z}(p^\infty)$ admits a neither maximally nor 
minimally almost-periodic Hausdorff group topology $\tau$
(Theorem~\ref{thm:non-map-MAP}). Thus,
$\mathbf{n}(\mathbb{Z}(p^\infty),\tau)$ is finite, and 
$\mathbf{n}(\mathbf{n}(\mathbb{Z}(p^\infty),\tau)) \subsetneq
\mathbf{n}(\mathbb{Z}(p^\infty),\tau)$. In particular, 
Problem~\ref{prob:nnG} is meaningful.


\section{{\itshape T}-sequences in abelian groups}

\label{sect:abelian}

In this section, $A$ is an abelian group and $\underline a=\{a_k\}\subseteq A$ is 
a sequence in $A$. In what follows, we provide several sufficient conditions for 
$\{a_k\}$  to be a $T$-sequence in $A$. For $l,m\in \mathbb{N}$, one puts
\begin{equation}
A(l,m)_{\underline a}=\{ m_1 a_{k_1} + \cdots + m_h a_{k_h} \mid
m\leq k_1<\cdots <k_h,m_i\in\mathbb{Z}\backslash\{0\},\sum |m_i| \leq l\}.
\end{equation}
The Zelenyuk-Protasov criterion for $T$-sequences states:

\begin{ftheorem}[{\cite[2.1.4]{PZMono}, \cite[Theorem~2]{ZelProt}}] \label{thm:ZelProt}
A sequence $\{a_k\}$ in an abelian group $A$ is a $T$-sequence if and 
only if for every $l\in\mathbb{N}$ and
$g\neq 0$, there exists $m \in \mathbb{N}$ such that
$g \not\in A(l,m)_{\underline a}$.
\end{ftheorem}

Put $A[n]=\{a\in A: na=0\}$ for every $n \in \mathbb{N}$. One says 
that $A$ is {\it almost torsion-free} if $A[n]$ is finite for every 
$n \in \mathbb{N}$ (cf. \cite{TkaYasch}). 


\begin{theorem} \label{thm:T:ord}
Let $A$ be an abelian group, and let $\{a_k\}\subseteq A$ be a 
sequence such that $t_k:=o(a_k)$ is finite for every $k\in \mathbb{N}$. 
Consider the following statements:

\begin{list}{{\rm (\roman{enumi})}}
{\usecounter{enumi}\setlength{\labelwidth}{25pt}\setlength{\topsep}{-10pt}
\setlength{\itemsep}{-0pt} \setlength{\leftmargin}{25pt}}

\item
\begin{equation} \label{eq:1-k}
\lim\limits_{k\rightarrow \infty}
\frac{t_k}{\gcd (t_k,\operatorname{lcm}(t_1,\cdots,t_{k-1}))}
=\infty.
\end{equation}

\item
For every $l\in \mathbb{N}$,
\begin{equation}
\lim\limits_{m\rightarrow \infty}
\inf \left\{ \max\limits_{1 \leq i \leq l} \frac{t_{k_i}}
{\gcd(t_{k_i},\operatorname{lcm}(t_{k_1},\ldots,t_{k_{i-1}},
t_{k_{i+1}},\ldots,t_{k_l}))}
\biggm|
\underline k \in \mathbb{N}^l_{m<}\right\} = \infty,
\end{equation}
where $\mathbb{N}^l_{m<}=\{\underline k = (k_1,\ldots,k_l)\in 
\mathbb{N}^l \mid m\leq k_1<\cdots <k_l\}$.

\item
For every $l\in \mathbb{N}$,
\begin{equation} \label{eq:ord}
\lim\limits_{m \rightarrow \infty}
\inf \{\max\limits_{1 \leq i \leq l} o(a_{k_i} + A_{\underline k^i}) \mid
\underline k \in \mathbb{N}^l_{m<} \} = \infty,
\end{equation}
where $A_{\underline k^i}=\langle a_{k_1},\ldots,a_{k_{i-1}},a_{k_{i+1}},
\ldots, a_{k_l}\rangle$.

\item
For every $l,n \in \mathbb{N}$, there exists $m_0\in \mathbb{N}$ such
that $A[n] \cap A(l,m)_{\underline a}=\{0\}$ for every $m \geq m_0$.

\item
$\{a_k\}$ is a $T$-sequence.
\end{list}
\vspace{10pt}
One has {\rm (i)} $\Rightarrow$ {\rm (ii)} $\Rightarrow$ {\rm (iii)} $\Rightarrow$ 
{\rm (iv)} $\Rightarrow$  {\rm (v)}, and if $A$
is almost torsion-free, then {\rm (v)} $\Rightarrow$
{\rm (iv)}.
\end{theorem}

\begin{proof}
(i) $\Rightarrow$  (ii) is obvious.

(ii) $\Rightarrow$ (iii):
Clearly, the order of $a_{k_i}$ in
$(\langle a_{k_i}\rangle + A_{\underline k^i})/A_{\underline k^i}$ 
is equal to its order modulo
$\langle a_{k_i} \rangle \cap A_{\underline k^i}$,\linebreak and 
$|\langle a_{k_i} \rangle \cap A_{\underline k^i}|$ divides both
$t_{k_i}$ and $\exp(A_{\underline k^i})$. The exponent
$\exp(A_{\underline k^i})$, in turn, divides\linebreak
$d=\operatorname{lcm}(t_{k_1},\ldots,t_{k_{i-1}}, t_{k_{i+1}},\ldots,t_{k_l})$,
because  $A_{\underline k^i}$ is generated by elements of orders
$t_{k_1},\ldots,t_{k_{i-1}},\linebreak[2] t_{k_{i+1}},\ldots,t_{k_l}$.
Therefore, $|\langle a_{k_i} \rangle \cap A_{\underline k^i}|$
divides their greatest common divisor of $t_k$ and $d$. Hence,
\begin{equation}
\frac{t_{k_i}}
{\gcd(t_{k_i},\operatorname{lcm}(t_{k_1},\ldots,t_{k_{i-1}},
t_{k_{i+1}},\ldots,t_{k_l}))} \biggm|   \frac{|\langle a_{k_i} \rangle |}
{|\langle a_{k_i} \rangle \cap A_{\underline k^i}|}
=o(a_{k_i} + A_{\underline k^i}).
\end{equation}

(iii) $\Rightarrow$ (iv): Given $l,n\in\mathbb{N}$, let $m_0\in\mathbb{N}$ be
such that $nl<\max\limits_{1\leq i\leq h} o(a_{k_i} + A_{\underline k^i})$
for every $1 \leq h \leq l$ and every $\underline k \in \mathbb{N}^h_{m_0<}$.
(By (\ref{eq:ord}), such $m_0$ exists.) Let 
$g=m_1 a_{k_1} + \cdots + m_h a_{k_h} \in A(l,m)_{\underline a}$ be a non-zero 
element, where $m_0 \leq m\leq k_1<\cdots <k_h$, $m_i\in\mathbb{Z}\backslash\{0\}$, 
and $\sum |m_i| \leq l$. It follows from the last two conditions that
$h \leq l$. So, there exists $1 \leq i \leq h$ such that
$nl < o(a_{k_i} + A_{\underline k^i})$, and thus
$n <  o(m_i a_{k_i} + A_{\underline k^i})$. To complete the proof, note
that $g \in m_i a_{k_i} + A_{\underline k^i}$, and therefore
$o(m_i a_{k_i} + A_{\underline k^i}) \mid o(g)$. Hence, 
$n < o(g)$, and so $g \not \in A[n]$, as desired.

(iv) $\Rightarrow$ (v):
Let $g \in A$ be a non-zero element. If the order of $g$ is infinite, then
$g \not\in A(l,1)_{\underline a}$ for every $l \in \mathbb{N}$, and so suppose
that $n:=o(g)$ is finite. By (iii), for every $l \in \mathbb{N}$ there exists
$m_0(l)$ such that $A[n] \cap A(l,m_0(l))_{\underline a}=\{0\}$.
In particular, $g \not\in A(l,m_0(l))_{\underline a}$ for every $l$.

(v) $\Rightarrow$ (iv):
Given $l,n\in \mathbb{N}$, and suppose that
$A[n]=\{0,g_1,\ldots,g_j\}$ is finite. For each $g_i$, pick
$m_i(l)\in\mathbb{N}$ such that $g_i \not\in A(l,m_i(l))_{\underline a}$, 
and put $m_0(l) = \max m_i(l)$. Clearly, one has 
$A[n] \cap A(l,m)_{\underline a} = \{0\}$
for every $m\geq m_0(l)$, as desired.
\end{proof}

\begin{remark}
In Theorem~\ref{thm:T:ord}, (iv) does not imply (iii). Indeed, although (iii)
fails for the sequence $\{e_n\}$ from Example~\ref{ex:p-n}, it is a
$T$-sequence in $\mathbb{Z}(p^\infty)$.
\end{remark}


\begin{corollary} \label{cor:div:T}
Let $A$ be an abelian group, and let $\{a_k\}\subseteq A$ be a 
sequence such that $t_k:=o(a_k)$ is finite for every $k\in \mathbb{N}$.

\begin{list}{{\rm (\alph{enumi})}}
{\usecounter{enumi}\setlength{\labelwidth}{25pt}\setlength{\topsep}{-10pt}
\setlength{\itemsep}{-4pt} \setlength{\leftmargin}{25pt}}

\item
If the $t_k$ are pairwise coprime, then $\{a_k\}$ is a $T$-sequence.

\item
If $t_k \mid t_{k+1}$ and $\lim\limits_{k \rightarrow \infty}
\frac {t_{k+1}}{t_k} = \infty$, then $\{a_k\}$ is a $T$-sequence.
\end{list}
\end{corollary}

\begin{proof}
If the $t_k$ are pairwise coprime, then the expression in (\ref{eq:1-k})
is equal to $t_k$ and $t_k \rightarrow \infty$. If $t_k \mid t_{k+1}$,
then the expression in (\ref{eq:1-k}) is precisely $\frac {t_{k+1}}{t_k}$.
In both cases, the statement follows from Theorem~\ref{thm:T:ord}(i).
\end{proof}

\section{Direct sums of finite abelian groups}

\label{sect:sum}

In this section, we provide a partial answer to Problem~\ref{prob:nG}:

\begin{theorem} \label{thm:sum:main}
Let $A=\bigoplus\limits_{\alpha \in I} F_\alpha$ be the direct sum
of an infinite family  $\{F_\alpha\}$ of non-trivial finite abelian
groups.
There exists a $T$-sequence $\{d_k\}$ in $A$ such that 
$A\{d_k\}$ is almost maximally-almost periodic.
\end{theorem}

\begin{remark}
In the setting of Theorem~\ref{thm:sum:main}, $A$ is obviously not 
algebraically  directly indecomposable. Furthermore, $A\{d_k\}$ need not 
be topologically directly indecomposable either: Consider the group
$B=\mathbb{Z}/2\mathbb{Z} \oplus \bigoplus\limits_{n=1}^\infty 
\mathbb{Z}/3\mathbb{Z}$, and let $\tau$ be a Hausdorff group topology on 
$B$. The subgroup $B_2 =\bigoplus\limits_{n=1}^\infty \mathbb{Z}/3\mathbb{Z}$  
is closed in $\tau$, because it is the kernel 
of the continuous group homomorphism $x \mapsto 3 x$. Thus, 
$(B,\tau)$ decomposes into a topological direct product 
of $B_1=\mathbb{Z}/2\mathbb{Z}$ and $B_2$ (where $B_1$ and $B_2$ 
are equipped with the subgroup topology). 
This also shows that $B_1 \cap\mathbf{n}(B,\tau)=\{0\}$, 
because $(B,\tau) \rightarrow B/B_2
\cong \mathbb{Z}/2\mathbb{Z}$ is continuous. In particular, not every 
finite subgroup of an abelian group $A$ is of the form 
$\mathbf{n}(A,\tau)$, where $\tau$ is a Hausdorff group topology on $A$.
\end{remark}

In order to prove Theorem~\ref{thm:sum:main}, we the following result:

\begin{proposition} \label{prop:sum:x}
Let $A=\bigoplus\limits_{i=1}^\infty C_i$ be a direct sum of cyclic groups 
of order $n_i=|C_i|>1$, and suppose that

\begin{list}{{\rm (\alph{enumi})}}
{\usecounter{enumi}\setlength{\labelwidth}{25pt}\setlength{\topsep}{-10pt}
\setlength{\itemsep}{-4pt} \setlength{\leftmargin}{25pt}}

\item
$n_i = n_{i+1}$ for every $i$.

\item
$n_i < n_{i+1}$ for every $i$, or

\vspace{11pt}
\end{list}
Then, for every $x \in A$, there is a $T$-sequence $\{d_k\}$ such 
that $\mathbf{n}(A\{d_k\})$ is finite and contains $x$.
\end{proposition}

\begin{proof}
The construction below is a modification of \cite[Example~5]{ZelProt} 
and \cite[2.6.2]{PZMono}. The sequence $d_k$ is constructed identically
in both (a) and (b), and the two are distinguished only in the proof of
$\{d_k\}$ being a $T$-sequence.

For each $i$, pick a generator $g_i$ in $C_i$. Each $y \in A$ can be 
written as $y=\sum \alpha_i g_i \in A$, and the $\alpha_i$ are unique modulo 
$n_i$. We set
$\Lambda(y) =\{ i \in\mathbb{N} \mid \alpha_i \not \equiv 0 \mod n_i\}$ 
and $\lambda(y)=|\Lambda(y)|$.  Put $i_0 = \max{\Lambda(x)}$. 
We define two sequences:
\begin{align}
a_k\colon &&& g_{i_0+1}, 2g_{i_0+1},\ldots, (n_{i_0+1}-1)g_{i_0+1},
g_{i_0+2},2g_{i_0+2},\ldots, (n_{i_0+2}-1)g_{i_0+2}, \ldots \\
b_k \colon &&& -x+ g_{i_0+1}, -x + g_{i_0+2}+g_{i_0+3},
-x + g_{i_0+4}+g_{i_0+5}+g_{i_0+6}, \ldots 
\end{align}
Let $\chi\colon A \rightarrow \mathbb{T}$ be a character of $A$. 
If $\chi$ is zero on all but finitely many of the $C_i$ and $\chi(x)=0$, 
then $\chi(a_k)=0$ and $\chi(b_k)=0$ for $k$ large enough, and so
$\chi(a_k)\longrightarrow 0$ and $\chi(b_k)\longrightarrow 0$. Conversely, 
suppose that $\chi(a_k)\longrightarrow 0$ and $\chi(b_k)\longrightarrow 0$.
Then there is $k_0\in \mathbb{N}$ such that 
$\chi(a_k) \subseteq (-\frac 1 3,\frac 1 3)$  for every $k>k_0$. Thus, 
there is $j_0\in \mathbb{N}$ such that $\chi(C_j) \subseteq 
(-\frac 1 3,\frac 1 3)$ for every $j > j_0$. Since the only subgroup 
contained in $(-\frac 1 3,\frac 1 3)$ is $\{0\}$, this means that $\chi$ 
is zero on all but finitely many of the $C_i$. Therefore, 
$\chi(b_k)=-\chi(x)$  for $k$ large enough, and hence $\chi(x)=0$.

The foregoing argument shows that if $\{d_k\}$ is any combination of the 
sequences $\{a_k\}$ and $\{b_k\}$ without repetitions (such as $a_1, b_1, 
a_2, b_2, \ldots$)
and if $\{d_k\}$ is a $T$-sequence, then $\chi$ is a continuous character
of $A\{d_k\}$ if and only if $\chi$ is zero on all but finitely many of 
the $C_i$ and $\chi(x)=0$. Thus, $x \in \mathbf{n}(A\{d_k\})$,
and the character $\chi_j\colon A \rightarrow \mathbb{T}$ defined by
$\chi_j(\sum\alpha_i g_i)=\frac 1 {n_j} \alpha_i$ is continuous on
$A\{d_k\}$ for  every $j > i_0$. Therefore, 
$\mathbf{n}(A\{d_k\})\subseteq C_1 \oplus \cdots \oplus C_{i_0}$, 
and hence $\mathbf{n}(A\{d_k\})$ is finite, as desired.

We show that $d_k$ is a $T$-sequence. First, observe that for every 
$l \in \mathbb{N}$ and every $j > i_0$ there exists 
$m\in \mathbb{N}$ such that 
\begin{align}
A(l, m)_{\underline d} & \subseteq \langle x \rangle \oplus
\bigoplus\limits_{i=j}^\infty C_i. \\
\intertext{Thus, for every $l \in \mathbb{N}$,}
\bigcap\limits_{m=1}^\infty A(l,m)_{\underline d} & \subseteq
\bigcap\limits_{j>i_0} \left(\langle x \rangle \oplus
\bigoplus\limits_{i=j}^\infty C_i\right) = \langle x \rangle.
\end{align}
Therefore, the condition of the Zelenyuk-Protasov criterion 
(Theorem~\ref{thm:ZelProt}) holds for every 
$y \not\in \langle x\rangle$, and it remains to show it for non-zero
elements of $\langle x\rangle$. Let $l \in \mathbb{N}$, and for
the time being assume only that $m >i_0$. 
If $\alpha x \in A(l,m)_{\underline d}$, then
\begin{align}
\alpha x & = (m_1 d_{k_1}+\cdots m_{h_1} d_{k_{h_1}}) + 
(m_{h_1+1} d_{k_{h_1+1}} + \cdots + m_h d_{k_h}), \\
\intertext{where $\sum |m_i| \leq l$, $m_i \neq 0$, $k_i \geq m$,
$d_{k_i}$ is a member of $\{a_k\}$ for $1 \leq i \leq h_1$,
and of  $\{b_k\}$ for $h_1+1\leq i \leq h$. (Here, we only assume that the 
$k_i$ are distinct, but they need not be increasing.) Thus,}
\alpha x & = (m_1 d_{k_1}+\cdots m_{h_1} d_{k_{h_1}}) +
(m_{h_1+1} (d_{k_{h_1+1}}+x) + \cdots + m_h (d_{k_h}+x)) -
\sum\limits_{i={h_1+1}}^h m_i x
\label{eq:sum:zero}
\end{align}
Since $k_i \geq m > i_0$, the first and the second expression on the right 
side belong to $\bigoplus\limits_{j=i_0+1}^\infty C_j$, while the left 
side and the third summand on the right belong to $\langle x \rangle$.
Therefore, 
\begin{align}
(m_1 d_{k_1}+\cdots m_{h_1} d_{k_{h_1}}) +
(m_{h_1+1} (d_{k_{h_1+1}}+x)+ \cdots + m_h (d_{k_h}+x)) = 0.
\label{eq:sum:zero2}
\end{align}
The sets  $\Lambda(m_{h_1+j} (d_{k_{h_1+j}}+x))$ are  disjoint, because 
the $k_i$ are distinct, and so 
\begin{align}
\lambda(m_{h_1+1} (d_{k_{h_1+1}}+x)+ \cdots + m_h (d_{k_h}+x)) =
\sum\limits_{j=1}^{h-h_1} \lambda(m_{h_1+j} (d_{k_{h_1+j}}+x)).
\label{eq:sum:lambda}
\end{align}
Since $\lambda(a_k)=1$ for every $k$, one has
$\lambda(m_1 d_{k_1}+\cdots m_{h_1} d_{k_{h_1}}) \leq h_1 \leq l$, 
and  hence, by 
(\ref{eq:sum:zero2}) and (\ref{eq:sum:lambda}),
\begin{align} \label{eq:lambda:contr}
\lambda(m_{h_1+j} (d_{k_{h_1+j}}+x))\leq l 
\end{align}
for every $1 \leq j\leq h - h_1$.

(a) Pick $m > i_0$ such that $\lambda(d_k +x) \geq l+1$ for every 
$k \geq m$ such that $d_k$ is a member of $\{b_k\}$. Then each 
$d_{k_{h_1+j}}+x$ is a sum of at least $l+1$ distinct base elements $g_i$ 
of order
$n=n_i$, and so $\lambda(m_{h_1+j} (d_{k_{h_1+j}}+x))\leq l$ 
implies $m_{h_1+j} (d_{k_{h_1+j}}+x)=0$. Therefore, $n \mid m_{h_1 +j}$, 
and in particular, $n \mid \sum\limits_{i={h_1+1}}^h m_i$. Hence,
by (\ref{eq:sum:zero}), $\alpha x = 0$, as desired.

(b) Pick $m$ as in (a), but with the additional condition that 
$n_i > l$ for every $i \in \Lambda(d_k+x)$  and for every $k \geq m$ such that 
$d_k$ is a member of $\{b_k\}$. This is possible because the 
$\Lambda(b_k+x)$ are disjoint and $\{n_i\}$ is increasing.
Each $d_{k_{h_1+j}}+x$ is a sum of at least $l+1$ distinct base elements, 
but of different orders. Thus, $m_{h_1+j} g_i \neq 0$ for every
$i \in  \Lambda(d_{k_{h_1+j}}+x)$ and $1 \leq j\leq h-h_1$,
because $|m_{h_1+j}| \leq l < n_i$ (by 
our assumption), and $m_{h_1+j} \neq 0$. Therefore,
\begin{equation}
\lambda(m_{h_1+j}(d_{k_{h_1+j}}+x)) = 
\sum\limits_{i \in   \Lambda(d_{k_{h_1+j}}+x)} 
\lambda(m_{h_1+j} g_i ) \geq l+1,
\end{equation}
contrary to (\ref{eq:lambda:contr}). Hence, $h=h_1$, and $\alpha x=0$, as 
desired. 
\end{proof}

\begin{proof}[Proof of Theorem~\ref{thm:sum:main}.]
Since each $F_\alpha$ decomposes into a direct sum of cyclic subgroups, we 
may assume that the $F_\alpha$ are cyclic from the outset. The set
$\{ |F_\alpha| \mid \alpha \in I\}$ is either bounded or contains an 
increasing sequence.  Thus, there is a countable subset $J \subseteq I$ 
such that $A_1=\bigoplus\limits_{\alpha \in J} F_\alpha$ satisfies one of 
the conditions of Proposition~\ref{prop:sum:x}. Pick $\gamma \in J$, and
let $g$ be a generator of $F_\gamma$. Then, by 
Proposition~\ref{prop:sum:x}, there is a $T$-sequence $\{d_k\}$ in $A_1$ such 
that $\mathbf{n}(A_1\{d_k\})$ is finite and contains $g$. For 
$A_2 = \bigoplus\limits_{\alpha \in I \backslash J} F_\alpha$, one has
$A\{d_k\}=A_1\{d_k\} \times A_2$, where $A_2$ is equipped with the 
discrete topology. Therefore,
$\mathbf{n}(A\{d_k\})=\mathbf{n}(A_1\{d_k\}) \times \mathbf{n}(A_2)=
\mathbf{n}(A_1\{d_k\})$ is non-trivial and finite, as desired.
\end{proof}

\section{Special sequences in the Pr\"ufer groups}

\label{sect:Pruefer}

In this section, we present method for constructing an algebraically 
directly indecomposable Hausdorff
abelian topological groups $A$ such that $\mathbf{n}(A)$ is
non-trivial and finite. An implicit yet rather thick hint for the 
construction 
of a group with these properties appears in \cite{DikMilTon}, in the 
proof of  Corollary~4.9 and the paragraph following it. It was 
Dikranjan who pointed out to the author that considering a suitable 
$T$-sequence in a Pr\"ufer group (and the maximal Hausdorff
group topology thus obtained) would lead to the desired example.
Pr\"ufer groups are distinguished by the property of 
having only finite proper subgroups, which implies that these subgroups 
are closed in any Hausdorff group topology. This property makes
Pr\"ufer groups  particularly suitable candidates for the 
aforesaid construction, because if 
$\mathbf{n}(\mathbb{Z}(p^\infty)\{d_n\})$ is a non-trivial proper 
subgroup, then it must be finite. Therefore, in this section, we study 
certain sequences in the Pr\"ufer groups  $\mathbb{Z}(p^\infty)$, and 
construct a $T$-sequence $\{d_n\}$ such that $\mathbb{Z}(p^\infty)\{d_n\}$
is neither maximally nor minimally almost-periodic. A second important
property that Pr\"ufer groups, being $p$-groups, have is that
for every $a,b \in \mathbb{Z}(p^\infty)$,
\begin{equation} \label{eq:p:ord}
o(a) \neq o(b) \Longrightarrow o(a+b)=\max\{o(a),o(b)\}.
\end{equation}
One says that a group $G$ is {\em potentially compact} if
for every ultrafilter $\mathcal{U}$ on $G$ there is $x\in G$
such that $\mathcal{U} -x$ is a $T$-filter, that is, 
$\mathcal{U}\stackrel\tau\longrightarrow x$
in some group topology $\tau$ (cf. \cite{Prota}, \cite{ZelProt2}).
A third noteworthy property is that Pr\"ufer groups  
are {\em not} potentially compact, because they are 
divisible torsion groups (cf. \cite[Thm.~6]{ZelProt2}).

\bigskip

Recall that if $A$ is a subgroup 
of an abelian Hausdorff topological group $S$, and $\{a_k\}\subseteq A$ 
is so that $a_k \longrightarrow b$ in $S$, where 
$\langle b\rangle \cap A = \{0\}$, then $\{a_k\}$ is a $T$-sequence
in $A$ (cf.~\cite[2.1.5]{PZMono}, \cite[Thm.~3]{ZelProt}). 
The setting of this result is so
that the sequence $a_n$ converges to an {\em external} element (namely, 
$b \not \in A$) in some group topology. In contrast, in this section,  
we investigate sequences in $\mathbb{Z}(p^\infty)$ that converge to a 
non-zero (internal) {\em element of} $\,\mathbb{Z}(p^\infty)$ 
in the ``usual" topology, that is, the one inherited from 
$\mathbb{Q}/\mathbb{Z}$.

\bigskip

We start off with an immediate consequence of Corollary~\ref{cor:div:T}.

\begin{lemma} \label{lemma:Pruefer:triv}
Let $\{a_k\}$ be a sequence in $\mathbb{Z}(p^\infty)$, and suppose that
$o(a_k)=p^{n_k}$. If $n_{k+1}-n_k \longrightarrow \infty$, then
$\{a_k\}$ is a $T$-sequence.
\end{lemma}

\begin{proof}
One has $\frac{o(a_{k+1})}{o(a_k)} = \frac{p^{n_{k+1}}}{p^{n_k}} = 
p^{n_{k+1} -n_k}$. Therefore, the statement follows from
Corollary~\ref{cor:div:T}(b).
\end{proof}

Example~\ref{ex:p-n} shows that the converse of Lemma~\ref{lemma:Pruefer:triv}
does not hold in general. Nevertheless, for some special sequences,
the condition of $n_{k+1}-n_k \longrightarrow \infty$ turns out to be
necessary for being a $T$-sequence, as Theorem~\ref{thm:Pruefer:spec} below
reveals.

We proceed by introducing some terminology.
A topological group $A$ is {\em precompact} if for every neighborhood $U$
of zero there is a finite subset $F \subseteq A$ such that
$A=F+U$. Following \cite{BarDikMilWeb}, we say that
a sequence $\{a_n\}$ on a group $G$ is a 
{\em $TB$-sequence} if there is a  precompact Hausdorff  group topology
$\tau$ on $G$ such that $a_n \stackrel\tau\longrightarrow 0$. 

It is easy to see that $A$ is precompact if and only if it carries the initial
topology induced by its group of continuous characters. Thus,
if $\{a_k\}$ is a sequence in an abelian group $A$, then 
by the universal property of $A\{a_k\}$, a character
$\chi \colon A \rightarrow \mathbb{T}$ is continuous on $A\{a_k\}$
if and only if $\chi(a_k)\longrightarrow 0$. Therefore, 
for $H= \{\chi \in \hom_\mathbb{Z}(A,\mathbb{T}) \mid
\chi(a_k)\longrightarrow 0 \}$, 
the closure of $\{0\}$ in the initial topology induced by $H$ is
$\mathbf{n}(A\{a_k\}) = \bigcap\limits_{\chi \in H} \ker \chi$.
Hence, $\{a_k\}$ is a $TB$-sequence if and only if
$H$ separates the points of $A$. (Observe that
$H=\widehat{A\{a_k\}}_d$.)

\begin{theorem} \label{thm:Pruefer:spec}
Let $x \in \mathbb{Z}(p^\infty)$ be a non-zero element, 
$\{n_k\} \subseteq\mathbb{N}$ an increasing sequence of positive integers, 
and set
\begin{equation}
a_k=-x + e_{n_k} = -x + \frac 1 {p^{n_k}} \in \mathbb{Z}(p^\infty).
\end{equation}

\begin{list}{{\rm (\alph{enumi})}}
{\usecounter{enumi}\setlength{\labelwidth}{25pt}\setlength{\topsep}{-10pt}
\setlength{\itemsep}{-4pt} \setlength{\leftmargin}{25pt}}

\item
$\{a_k\}$ is a $T$-sequence if and only if
$n_{k+1}-n_k\longrightarrow\infty$.

\item
$\{a_k\}$ is a $T$-sequence if and only if it is a $TB$-sequence. In this
case, $\mathbb{Z}(p^\infty)\{a_k\}$ is maximally  almost-periodic, and it
has $\mathfrak{c}$ many faithful characters {\rm (}in  particular,
$|\widehat{\mathbb{Z}(p^\infty)\{a_k\}}|=\mathfrak{c}${\rm )}.
\end{list}
\end{theorem}

Since every $TB$-sequence is a $T$-sequence, but the converse is not true
in general, (b) of the Theorem~\ref{thm:Pruefer:spec} is a non-trivial
result. Its proof, however, requires a technical lemma. Note that the
group of all characters of $\mathbb{Z}(p^\infty)$ is isomorphic
to the group $\mathbb{Z}_p$ of the $p$-adic integers. In other words, 
$\mathbb{Z}_p=\hom_\mathbb{Z}(\mathbb{Z}(p^\infty),\mathbb{T})$.

\begin{lemma} \label{lemma:r_k}
Suppose that $n_{k+1}-n_k \longrightarrow \infty$. For 
$\chi=\sum\limits_{n=0}^\infty\alpha_n p^n\in\mathbb{Z}_p$
{\rm (}$0 \leq \alpha_n \leq p - 1${\rm )} and $\gamma \in (0,1)$,
$\chi(e_{n_k}) \longrightarrow \gamma$ if and only if
\begin{equation}
r_k :=
\frac{\sum\limits_{l=n_k}^{n_{k+1}-1}\alpha_l p^{l-n_k}}{p^{n_{k+1}-n_k}}
\longrightarrow \gamma.
\end{equation}
\end{lemma}

\begin{proof}
One has
\begin{align}
\chi(e_{n_{k+1}}) & = \frac{\sum\limits_{l=0}^{n_{k+1}-1} \alpha_l p^l}
{p^{n_{k+1}}} =
\frac{\sum\limits_{l=0}^{n_{k}-1} \alpha_l p^l} {p^{n_{k+1}}} +
\frac{\sum\limits_{l=n_k}^{n_{k+1}-1} \alpha_l p^l}{p^{n_{k+1}}}=
\frac{\chi(e_{n_{k}})} {p^{n_{k+1}-n_k}} + r_k,
\end{align}
and thus $\lim\limits_{k \rightarrow\infty} \chi(e_{n_k}) =
\lim\limits_{k \rightarrow\infty} r_k$ in $\mathbb{T}$
(by the equality of limits we mean that one exists if and only
if the other does, and in that case they are equal),
because $n_{k+1}-n_k \longrightarrow \infty$. Since $\gamma \neq 0$,
small enough neighborhoods of $\gamma$ in $\mathbb{T}$ and $(0,1)$ are the
same, and therefore $\lim\limits_{k \rightarrow\infty} \chi(e_{n_k})
= \lim\limits_{k \rightarrow\infty} r_k$ in $(0,1)$.
\end{proof}

\begin{proof}[Proof of Theorem~\ref{thm:Pruefer:spec}]
(a) Since $\{n_k\}$ is increasing, one has 
$n_k \longrightarrow \infty$. Thus, $p^{n_k} > o(x)$ for $k$ large enough,
and so $o(a_k)=p^{n_k}$ except for maybe a finite number of $k$
(by (\ref{eq:p:ord})). Hence, $a_k$ is a $T$-sequence by 
Lemma~\ref{lemma:Pruefer:triv}.

Conversely, let $p^{n_0}=o(x)$, and assume that
$n_{k+1}-n_k\not\longrightarrow\infty$. Then  $o(p^{n_0-1}x)=p$,
\begin{equation}
p^{n_0-1} a_k = -p^{n_0-1}x + e_{n_k -n_0+1},
\end{equation}
and the differences $(n_{k+1}-n_0+1)-(n_k-n_0+1) = n_{k+1}-n_k
\not\longrightarrow \infty$. Thus, it suffices to show that $p^{n_0-1}a_k$
is not a $T$-sequence. Therefore,
without loss of generality, we may assume that $o(x)=p$ from the outset.
Since $n_{k+1}-n_k\not\longrightarrow\infty$, there exists a number $d$ 
and a subsequence $k_r$ of $k$ such that $n_{k_r +1} -n_{k_r} \leq d$ 
for every $r$. If $a_k \longrightarrow 0$ in a group topology $\tau$ on 
$\mathbb{Z}(p^\infty)$, then
in particular, $p a_{n_{k_r+1}} = e_{n_{k_r+1}-1} \longrightarrow 0$,
and so for every $1 \leq i \leq d$, $e_{n_{k_r+1} -i} \longrightarrow 0$.
Thus, the sequence $b_n$ defined as
\[
e_{n_{k_1+1} -d}, e_{n_{k_1+1} -d+1},\ldots,e_{n_{k_1+1}-1},
e_{n_{k_2+1} -d}, e_{n_{k_2+1} -d+1},\ldots,e_{n_{k_2+1}-1},\ldots
\]
also converges to $0$ in $\tau$. One has
$n_{k_r+1} -d \leq n_{k_r} \leq n_{k_r+1} -1$, and therefore
$e_{n_{k_r}}$ is a subsequence of $b_n$, and hence
$e_{n_{k_r}} \longrightarrow 0$ in $\tau$. Since
$a_{k_r}=- x + e_{n_{k_r}}$, this shows that $\tau$ is not Hausdorff.

(b) If $\{a_k\}$ is a $TB$-sequence, then clearly it is a $T$-sequence. 
Conversely, suppose that $\{a_k\}$ is a $T$-sequence. In order to 
show that $\{a_n\}$ is a $TB$-sequence, 
we find a faithful continuous character of $\mathbb{Z}(p^\infty)\{a_k\}$,
in other words, $\chi \in \mathbb{Z}_p$ such that
$\chi(a_k) \longrightarrow 0$ and $\ker \chi = \{0\}$. 
Let $p^{n_0}=o(x)$. 
For $\chi=\sum\limits_{n=0}^\infty\alpha_n p^n \in \mathbb{Z}_p$, 
a character of $\mathbb{Z}(p^\infty)$, if 
\begin{equation} \label{eq:cond:1}
\alpha_0=1,\alpha_1=\cdots=\alpha_{n_0 -1}=0,
\end{equation} 
then $\chi$ acts on the subgroup $\langle x \rangle$ as the identity, where
$\mathbb{Z}(p^\infty)$ is viewed as a subgroup of $\mathbb{T}$.
Thus,  $\chi(e_1)=e_1\neq 0$ and $\chi(x)=x \neq 0$,
and in  particular, $\chi$ is faithful. By Lemma~\ref{lemma:r_k}, 
$\chi(a_k)\longrightarrow 0$ (i.e., $\chi$ is continuous on 
$\mathbb{Z}(p^\infty)\{a_k\}$) if and only if
\begin{equation} \label{eq:cond:2}
r_k =
\frac{\sum\limits_{l=n_k}^{n_{k+1}-1}\alpha_l p^{l-n_k}}{p^{n_{k+1}-n_k}}
\longrightarrow x.
\end{equation}
Conditions (\ref{eq:cond:1}) and (\ref{eq:cond:2}) are satisfied 
(simultaneously) by continuum many elements in $\mathbb{Z}_p$, which 
completes the proof.
\end{proof}

We proceed by presenting the construction of a non-minimally almost-periodic 
non-maximally almost-periodic Hausdorff group topology on the group 
$\mathbb{Z}(p^\infty)$ for $p \neq 2$. Our technique makes substantial use 
of the assumption that $p\neq 2$; nevertheless, we conjecture that a similar
construction is available for $p=2$.

\begin{theorem} \label{thm:non-map-MAP}
Let $p$ be a prime number such that $p\neq 2$,
$x \in \mathbb{Z}(p^\infty)$ be a non-zero element with $p^{n_0}=o(x)$,
and put
\begin{equation}
b_n= -x +e_{n^3 -n^2} + \cdots + e_{n^3-2n} + e_{n^3-n} +e_{n^3}
=-x + \frac 1 {p^{n^3 -n^2}}
+ \cdots + \frac 1 {p^{n^3-2n}} + \frac 1 {p^{n^3-n}} +\frac 1 {p^{n^3}}.
\end{equation}
Consider the sequence $d_n$ defined as $b_1,e_1,b_2,e_2,b_3,e_3,\ldots$.
Then:
\begin{list}{{\rm (\alph{enumi})}}
{\usecounter{enumi}\setlength{\labelwidth}{25pt}\setlength{\topsep}{-10pt}
\setlength{\itemsep}{-4pt} \setlength{\leftmargin}{25pt}}

\item
$\{d_n\}$ is a $T$-sequence in $\mathbb{Z}(p^\infty)$;

\item
the underlying group of $\widehat{\mathbb{Z}(p^\infty)\{d_n\}}$ 
is $p^{n_0} \mathbb{Z}\subseteq \mathbb{Z}_p = 
\hom_\mathbb{Z}(\mathbb{Z}(p^\infty),\mathbb{T})$;

\item
$\mathbf{n}(\mathbb{Z}(p^\infty)\{d_n\})= \langle x\rangle$.
\end{list}

\vspace{11pt}

\noindent
In particular,  $\mathbb{Z}(p^\infty)\{d_n\}$ is neither maximally almost
periodic nor minimally almost-periodic, and
$\mathbf{n}(\mathbb{Z}(p^\infty)\{d_n\})$ is finite.
\end{theorem}

\begin{corollary} \label{cor:non-map-MAP}
Let $p$ be a prime number such that $p\neq 2$, and put
\begin{equation}
b_n= -e_1 +e_{n^3 -n^2} + \cdots + e_{n^3-2n} + e_{n^3-n} +e_{n^3}
=-\frac 1 p + \frac 1 {p^{n^3 -n^2}}
+ \cdots + \frac 1 {p^{n^3-2n}} + \frac 1 {p^{n^3-n}} +\frac 1 {p^{n^3}}.
\end{equation}
Consider the sequence $d_n$ defined as $b_1,e_1,b_2,e_2,b_3,e_3,\ldots$.
Then:
\begin{list}{{\rm (\alph{enumi})}}
{\usecounter{enumi}\setlength{\labelwidth}{25pt}\setlength{\topsep}{-10pt}
\setlength{\itemsep}{-4pt} \setlength{\leftmargin}{25pt}}

\item
$\{d_n\}$ is a $T$-sequence in $\mathbb{Z}(p^\infty)$;

\item
the underlying group of $\widehat{\mathbb{Z}(p^\infty)\{d_n\}}$ 
is $p \mathbb{Z}\subseteq \mathbb{Z}_p = 
\hom_\mathbb{Z}(\mathbb{Z}(p^\infty),\mathbb{T})$;

\item
$\mathbf{n}(\mathbb{Z}(p^\infty)\{d_n\})= \langle \frac 1 p\rangle$.
\end{list}

\vspace{11pt}

\noindent
In particular,  $\mathbb{Z}(p^\infty)\{d_n\}$ is neither maximally almost
periodic nor minimally almost-periodic, and
$\mathbf{n}(\mathbb{Z}(p^\infty)\{d_n\})$ is finite.
\end{corollary}

In order to prove Theorem~\ref{thm:non-map-MAP}, several auxiliary results of 
a technical nature are required. Until the end of this section, we assume 
that $p\neq 2$. Each element $y \in \mathbb{Z}(p^\infty)$ admits many 
representations of the form $y=\sum \sigma_n e_n$, where 
$\sigma_n \in \mathbb{Z}$ (only finitely many of the $\sigma_n$
are non-zero), and so we say that it  is the {\em canonical form} 
of $y$ if $|\sigma_n|\leq \frac{p-1}2$ for every $n \in \mathbb{N}$; 
in this case, we put $\Lambda(y)= \{ n \in \mathbb{N} \mid \sigma_n \neq 0\}$ 
and $\lambda(y)= | \Lambda(y)|$.

\begin{lemma} \label{lemma:Pruefer:canon}
Let $y =\sum \sigma_n e_n \in \mathbb{Z}(p^\infty)$. Then:

\begin{list}{{\rm (\alph{enumi})}}
{\usecounter{enumi}\setlength{\labelwidth}{25pt}\setlength{\topsep}{-10pt}
\setlength{\itemsep}{-4pt} \setlength{\leftmargin}{25pt}}

\item
$y$ admits a canonical form $y=\sum \sigma_n^\prime e_n$, and
$\sum |\sigma_n^\prime| \leq \sum |\sigma_n|$;

\item
the canonical form is unique, and so $\Lambda$ is well-defined.

\vspace{-4pt}
\hspace{-25pt} Furthermore, 

\vspace{-2pt}
\item
$\lambda(z) \leq l$ for every 
$z \in \mathbb{Z}(p^\infty)(l,1)_{\underline e }$ and $l \in \mathbb{N}$.
\end{list}
\end{lemma}

\begin{proof}
(a) Let $N$ be the largest index such that $\sigma_N \neq 0$. We proceed
by induction on $N$. If $N=1$, then $y=\sigma_1 c_1$. Thus, if
$\sigma_1=\sigma_1^\prime +  m p$ is a division with residue in $\mathbb{Z}$,
and $\sigma_1^\prime$ is chosen to have the smallest possible absolute value,
then  $|\sigma_1^\prime| \leq \frac{p-1} 2 $, and
\begin{equation}
y=(\sigma_1^\prime + p m)e_1 = \sigma_1^\prime e_1 +  m p e_1
= \sigma_1^\prime e_1.
\end{equation}
In particular, $|\sigma_1^\prime| \leq |\sigma_1|$.
Suppose now that the statement holds for all elements with representation
with maximal non-zero index less than $N$. If
$\sigma_N=\sigma_N^\prime +  k p$ is a division with residue in $\mathbb{Z}$,
and $\sigma_N^\prime$ is chosen to have the smallest possible absolute value,
then  $|\sigma_N^\prime| \leq \frac{p-1} 2 $, and
\begin{align}
\hspace{-2pt}
y-\!\sum\limits_{n=1}^{N-2}\sigma_n e_n  - (\sigma_{N-1}+k) e_{N-1} =
- k e_{N-1} +\sigma_N e_N
= - k e_{N-1} +(\sigma_N^\prime +  k p)e_N = \sigma_N^\prime e_N.
\end{align}
The element $z=\sum\limits_{n=1}^{N-2}\sigma_n e_n+(\sigma_{N-1}+k)e_{N-1}$
satisfies the inductive hypothesis, so
$z=\sum\limits_{n=1}^{N-1}\sigma_n^\prime e_n$,  where
$|\sigma_n^\prime|\leq \frac{p-1}2$ and
$\sum\limits_{n=1}^{N-1}|\sigma_n^\prime| \leq
\sum\limits_{n=1}^{N-1}|\sigma_n| +|k|$. Therefore,
$y=\sum \sigma_n^\prime c_n$, $|\sigma_n^\prime|\leq \frac{p-1}2$,
and $\sum |\sigma_n^\prime| \leq \sum |\sigma_n|$, because
$|\sigma_N^\prime|+|k| \leq |\sigma_N|$.

(b) Suppose that $\sum\sigma_n e_n = \sum \upsilon_n e_n$ are two
distinct canonical representations of the same element. Then
$\sum (\sigma_n -\upsilon_n) e_n = 0$, and $|\sigma_n - \upsilon_n|\leq p-1$.
Let $N$ be the largest index  such that $\sigma_N \neq \upsilon_N$.
(Since all coefficients are zero, except for a finite number of indices,
such $N$ exists.) This means that $0<|\sigma_N - \upsilon_N|\leq p-1$,
and $o((\sigma_N - \upsilon_N)e_N)=p^N$. Therefore, by
(\ref{eq:p:ord}), one has
$o(\sum (\sigma_n -\upsilon_n) e_n)=p^N$, because 
$o(\sum\limits_{n<N} (\sigma_n -\upsilon_n) e_n ) \leq p^{N-1}$. 
This is a contradiction, and therefore
$\sigma_n=\upsilon_n$ for every $n \in \mathbb{N}$.

(c) Let $z= \mu_1 e_{n_1} + \cdots +\mu_h e_{n_h}$,
where $\sum |\mu_i| \leq l$ and $n_1 < n_2 < \cdots < n_h$.
By (a), $z$ admits a canonical form $z=\sum \mu_n^\prime e_n$,
and  $\sum |\mu_n^\prime| \leq \sum |\mu_i| \leq l$. Therefore,
$\mu_n^\prime \neq 0$ only for at most $l$ many indices.
\end{proof}

\begin{lemma} \label{lemma:Pruefer:mcn}
Let $m\in \mathbb{Z}\backslash\{0\}$, and put $l=\lceil\log_p |m|\rceil$.
If $n > l$, then $\Lambda(m e_n)\subseteq  \{n-l,\ldots,n-1 ,n\}$ and
$1 \leq \lambda(m e_n)$.
\end{lemma}

\begin{proof}
It follows from $n >l$ that $p^n > |m|$, and so $me_n \neq 0$. Thus,
$1 \leq \lambda(m e_n)$. To show the first statement, expand
$m = \mu_0 + \mu_1 p + \cdots + \mu_l p^l$, where
$\mu_i \in \mathbb{Z}$ and $|\mu_i| \leq \frac {p-1} 2$. Then
\begin{equation}
m e_n = \mu_0 e_n + \mu_1 e_{n-1} + \cdots + \mu_l e_{n-l}
\end{equation}
is in canonical form, and therefore
$\Lambda(m e_n)\subseteq \{n-l,\ldots,n-1 ,n\}$, as desired.
\end{proof}

\begin{lemma} \label{lemma:Pruefer:y+z}
Let $y,z \in \mathbb{Z}(p^\infty)$ such that $\lambda(y) > \lambda(z)$,
and suppose that $\Lambda(y)=\{k_1,\ldots,k_g\}$
where $k_1 < \ldots < k_g$ and $g=\lambda(y)$.
Then $o(y+z) \geq p^{k_{g-\lambda(z)}}$.
\end{lemma}

\begin{proof}
Let $y=\sum \nu_n e_n$ and $z=\sum\mu_n e_n$ in canonical form.
Then $y+z = \sum (\nu_n + \mu_n) e_n$, and $|\nu_n+\mu_n|\leq p-1$.
Clearly, $o(y+z)=p^N$ for $N$ the largest index such
that $\nu_N +\mu_N \neq 0$. By the definition of $N$,
$\mu_n=-\nu_n$ for every $n>N$. In particular,
$\mu_{k_i}\neq 0$ for every $i$ such that $k_i > N$.
Thus, there are at most $\lambda(z)$ many $i$ such that
$k_i > N$, and therefore $N \geq k_{g-\lambda(z)}$.
\end{proof}

\begin{remark} \label{rem:Pruefer:lambda}
If $y_1,y_2 \in \mathbb{Z}(p^\infty)$ and
$\Lambda(y_1)\cap\Lambda(y_2)=\emptyset$,
then $\Lambda(y_1 + y_2) = \Lambda(y_1)\cup \Lambda(y_2)$ and
$\lambda(y_1 + y_2) = \lambda(y_1) + \lambda(y_2)$.
\end{remark}

\begin{proposition} \label{prop:Pruefer:lambda}
Let $y=\nu_1 e_{n_1}+\cdots \nu_f e_{n_f}$, where $n_1 < \cdots < n_f$
and $\nu_i\neq 0$.
Put $l_i=\lceil\log_p |\nu_i|\rceil$, and suppose that
$n_i < n_{i+1} -l_{i+1}$ for each $1\leq i\leq f$. Then:

\begin{list}{{\rm (\alph{enumi})}}
{\usecounter{enumi}\setlength{\labelwidth}{25pt}\setlength{\topsep}{-10pt}
\setlength{\itemsep}{-4pt} \setlength{\leftmargin}{25pt}}

\item
$f \leq \lambda(y)$;

\item
if $z \in \mathbb{Z}(p^\infty)$ is such that $\lambda(z) < \lambda(y)$, then
$o(y+z)\geq p^{n_{f-\lambda(z)}-l_{f-\lambda(z)}}$.
\end{list}
\end{proposition}

\begin{proof}
(a) By Lemma~\ref{lemma:Pruefer:mcn},
$\Lambda(\nu_i e_{n_i}) \subseteq \{n_i-l_i,\ldots,n_i\}$, and since
$n_{i-1} < n_i - l_i$, the sets $\Lambda(\nu_i e_{n_i})$ are pairwise
disjoint. Therefore, by Remark~\ref{rem:Pruefer:lambda},
$\lambda(y) = \lambda(\nu_1 e_{n_1})+\cdots + \lambda(\nu_f e_{n_f})\geq f$,
and
\begin{equation} \label{eq:prop:lambda}
\Lambda(y)\subseteq \bigcup\limits_{i=1}^f \{n_i-l_i,\ldots,n_i\} =
\{n_1-l_1,\ldots,n_1,\ldots, n_i-l_i,\ldots,n_i,\ldots,
n_f-l_f,\ldots,n_f\}.
\end{equation}

(b) By Lemma~\ref{lemma:Pruefer:y+z}, 
$o(y+z)\geq p^{k_{\lambda(y)-\lambda(z)}}$,
where $\Lambda(y)=\{k_1,\ldots,k_g\}$ (increasingly ordered).
Since $\Lambda(\nu_i e_{n_i})$ is
non-empty for each $i$, it follows from (\ref{eq:prop:lambda}) that
$k_{\lambda(y)-\lambda(z)} \geq n_{f-\lambda(z)}-l_{f-\lambda(z)}$.
\end{proof}

\begin{corollary} \label{cor:Pruefer:trick}
Let $l \in \mathbb{N}$, $z \in \mathbb{Z}(p^\infty)(l,1)_{\underline e }$, and
$y=e_{n_1}+\cdots+e_{n_f}$ such that $n_1 < \cdots < n_f$, $l < f$,
and $n_i  < n_{i+1} - l$. Then $o(\mu y + z) \geq p^{n_{f-l} -l}\geq p^{n_1-l}$
for every $\mu \in \mathbb{Z}$ such that $0<|\mu| \leq l$.
\end{corollary}

\begin{proof}
Since $|\mu| \leq l$, $\mu y = \mu e_{n_1}+\cdots+\mu e_{n_f}$ satisfies the 
conditions of Proposition~\ref{prop:Pruefer:lambda} (because 
$\log_p |\mu| \leq l$), and thus, $l<f \leq \lambda(\nu y)$. On the other hand, 
by Lemma~\ref{lemma:Pruefer:canon}(c), $\lambda(z) \leq l$, and therefore
$o(\nu y+z)\geq p^{n_{f-\lambda(z)}-l}\geq p^{n_{f-l}-l}$
pursuant to Proposition~\ref{prop:Pruefer:lambda}(b).
\end{proof}

\begin{proof}[Proof of Theorem~\ref{thm:non-map-MAP}.]
To shorten notations, put $A=\mathbb{Z}(p^\infty)$.

(a) In order to prove that $\{d_n\}$ is a $T$-sequence,
we show that (iv) of Theorem~\ref{thm:T:ord} holds.
For $n$ large enough, $o(b_n)=p^{n^3}$, and so
by Lemma~\ref{lemma:Pruefer:triv}, $\{b_n\}$ is a $T$-sequence; 
$\{e_n\}$ is evidently a $T$-sequence (cf.~Example~\ref{ex:p-n}). Thus,
by Theorem~\ref{thm:T:ord}, there exists $m_0$ such that
\begin{equation} \label{eq:map-MAP:1}
A[n] \cap A(l,m)_{\underline b} = A[n] \cap A(l,m)_{\underline e}=\{0\}
\end{equation}
for every $m \geq m_0 $ (because $A$ is almost torsion-free).
Without loss of generality, we may assume that
$m_0 > l+ n+n_0$. Observe that
\begin{equation}
A(l,2m)_{\underline d}\subseteq A(l,m)_{\underline b}\cup A(l,m)_{\underline e}
\cup (A(l,m)_{\underline b}\backslash\{0\}+
A(l,m)_{\underline e}\backslash\{0\}),
\end{equation}
and therefore it suffices to show that
$A(l,m)_{\underline b}\backslash\{0\}+A(l,m)_{\underline e}\backslash\{0\}$
contains no element of $A[n]$ for every $m\geq m_0$. Let
$z\in A(l,m)_{\underline e}\backslash\{0\}$ and
$w =m_1 b_{n_1}+\cdots + m_h b_{n_h}
\in A(l,m)_{\underline b}\backslash\{0\}$ where
$m\leq n_1 < \cdots < n_h$ and $0<\sum |m_i| \leq l$.
Put $y=e_{n_h^3 -n_h^2} + \cdots+ e_{n_h^3 -n_h} + e_{n_h^3}$. The number of 
summands in $y$ is $n_h+1$, and the differences between the indices of the 
terms is $n_h$. By the construction, $n_h \geq m \geq m_0 > l$. Thus, the 
conditions of Corollary~\ref{cor:Pruefer:trick} are satisfied, and since
$|m_h| \leq l$, we get $o(m_h y + z) \geq p^{n_h^3 - n_h^2 -l}> p^{(n_h-1)^3}$.
Therefore, $o(-m_h x) \neq o(m_h y + z)$ (because 
$o(-m_h x) \leq p^{n_0} \leq p^{m_0-1} \leq p^{(n_h-1)^3}$), and
\begin{equation} \label{eq:non-map-MAP:1}
o(m_h b_{n_h} + z) = o(-m_h x + m_h y + z) \stackrel{(\ref{eq:p:ord})}= 
\max\{o(-m_h x),o(m_h y + z)\} > p^{(n_h-1)^3}.
\end{equation}
One has
\begin{align}
o(w-m_h b_{n_h}) & \leq o(b_{n_{h-1}})=p^{n_{h-1}^3} \leq
p^{(n_h-1)^3}  \stackrel{(\ref{eq:non-map-MAP:1})}< o(m_h b_{n_h} + z), \\
\intertext{and hence}
o(w+z) & =o((w-m_h b_{n_h})+(z+m_h b_{n_h})) \\
& \hspace{-0.85pt} \stackrel{(\ref{eq:p:ord})} = 
\max \{o(w-m_h b_{n_h}),o(z+m_h b_{n_h})\} > p^{(n_h-1)^3}
> p^{(m_0-1)^3} >n.
\end{align}

(b) As noted earlier, a character
$\chi\in\hom_\mathbb{Z}(\mathbb{Z}(p^\infty),\mathbb{T})$ is continuous on 
$\mathbb{Z}(p^\infty)\{d_n\}$ if and only if $\chi(d_n) \longrightarrow 0$
(by the universal property)---in other words, $\chi(b_n)\longrightarrow 0$ 
and $\chi(e_n)\longrightarrow 0$. The latter is equivalent to $\chi$
having the form of $m \chi_1$, where $\chi_1$ is  the natural embedding of 
$\mathbb{Z}(p^\infty)$ into $\mathbb{T}$ and $m \in \mathbb{Z}$ 
(cf.~\cite[Example~6]{ZelProt}, \cite[3.3]{DikMilTon}). Since
\begin{equation}
0\leq \frac 1 {p^{n^3 -n^2}} + \cdots + \frac 1 {p^{n^3-2n}} +
\frac 1 {p^{n^3-n}} +\frac 1 {p^{n^3}} \leq \frac {n+1}{p^{n^3 -n^2}}
\longrightarrow 0,
\end{equation}
one has $\chi_1(b_n) \longrightarrow -x$, and consequently
$\chi(b_n)=m\chi_1(b_n)\longrightarrow 0$
if and only if $-mx=0$ (i.e., $x \in \ker \chi$). This means that
$\chi=m\chi_1$ if and only if
$o(x)=p^{n_0}\mid m$, as desired.

(c) We have already seen that $x \in \ker \chi$ for every continuous
character  of $\mathbb{Z}(p^\infty)\{d_k\}$. On the other hand,
$\mathbf{n}(\mathbb{Z}(p^\infty)\{d_k\}) \subseteq \ker p^{n_0} \chi_1
=\langle x \rangle$.
\end{proof}

\begin{remark} \label{rem:aMAP:ess}
A careful examination of the construction in Theorem~\ref{thm:non-map-MAP}
reveals that the only following properties of the sequence $\{b_n\}$ are
essential:

\begin{list}{{\rm (\arabic{enumi})}}
{\usecounter{enumi}\setlength{\labelwidth}{25pt}\setlength{\topsep}{-10pt}
\setlength{\itemsep}{-4pt} \setlength{\leftmargin}{25pt}}

\item
Growing number of summands in $b_n$---in other words, 
$\lambda(b_n)\longrightarrow \infty$;

\item
Growing gaps between the orders of summands in $b_n$ (in its
canonical form);

\item
$b_n \longrightarrow -x$ in the topology of inherited from
$\mathbb{T}$, where $x\in \mathbb{Z}(p^\infty)$ and
$x\neq 0$.

\vspace{6pt}
\end{list}
Condition (1) and (2) are needed in order to apply 
Corollary~\ref{cor:Pruefer:trick}, while (3) guarantees that 
$m\chi_1$ is continuous if and only if $o(x) \mid m$
(where $\chi_1$ is the natural embedding  of $\mathbb{Z}(p^\infty)$ into
$\mathbb{T}$).

\end{remark}

We conclude with a problem motivated by Theorem~\ref{thm:Pruefer:spec}
and Remark~\ref{rem:aMAP:ess}:

\begin{problem}
Is there a $T$-sequence $\{a_k\}$ in $\mathbb{Z}(p^\infty)$
with  bounded $\lambda(a_k)$ that is not a $TB$-sequence?
\end{problem}


\section*{Acknowledgments}

I am deeply indebted to Professor Dikran Dikranjan, who introduced me 
to this topic and spent so much time corresponding with me about it.
Without his attention and encouragement, I would have struggled to 
carry out this research.

I wish to thank Professor Horst Herrlich, Professor Salvador
Hern\'andez and Dr. Gavin Seal for the valuable discussions 
and helpful suggestions that were
of great assistance in writing this paper.  

I am grateful for the constructive comments of the anonymous referee that
led to a substantial improvement of the results and their presentation
in this paper.

\bibliography{notes,notes2,notes3}

\begin{samepage}

{\bigskip\bigskip\noindent 
Department of Mathematics and Statistics\\
Dalhousie University\\
Halifax, B3H 3J5, Nova Scotia\\
Canada

\nopagebreak
\bigskip\noindent{\em e-mail: lukacs@mathstat.yorku.ca} }
\end{samepage}

\end{document}